\documentclass[10pt,twoside,a4paper]{article}
\usepackage[italian,english]{babel}
\usepackage{amsmath,amsfonts,amssymb,amscd,lscape,amstext,amsthm}
\usepackage{tabularx}
\newtheorem{lemma}{\rm \indent LEMMA}[section]
\newtheorem{defin}{\rm \indent D\small EFINITION}[section]
\newtheorem{teore}[lemma]{\rm \indent T\small HEOREM}
\newtheorem{osse}[lemma]{\rm \indent Remark}
\newtheorem{corol}[lemma]{\rm \indent C\small OROLLARY}
\numberwithin{equation}{section}
\begin{document}
\title{\textbf{The local Calder\`{o}n problem and the determination
at the boundary of the conductivity}}
\author{\textsc{Giovanni Alessandrini}\thanks{Dipartimento di Matematica e Informatica, Universit\`{a} degli
Studi di Trieste, Via Valerio 12/b, 34127 Trieste,
alessang$@$units.it} and \textsc{Romina
Gaburro}\thanks{Department of Mathematics and Statistics,
University of Limerick, Castletroy, Limerick, Ireland,
romina.gaburro$@$ul.ie}}
\date{}
\maketitle \footnotesize{\textbf{Abstract.} We discuss the inverse
problem of determining the, possibly anisotropic, conductivity of
a body  $\Omega\subset\mathbb{R}^{n}$ when the so--called
Dirichlet-to-Neumann map is locally given on a non empty portion
$\Gamma$ of the boundary $\partial\Omega$. We extend results of
uniqueness and stability at the boundary, obtained by the same
authors in  \emph{SIAM J. Math. Anal.  33  (2001),  no. 1,
153--171}, where the Dirichlet-to-Neumann map was given on all of
$\partial\Omega$ instead. We also obtain a pointwise stability
result at the boundary among the class of conductivities which are
continuous at some point $y\in\Gamma$. Our arguments also apply
when the local Neumann-to-Dirichlet map is available.
\section{\normalsize{Introduction}.}
\normalsize In absence of internal sources, the electrostatic
potential $u$  in a conducting body, described by a domain
$\Omega\subset{\mathbb R}^n$, is governed  by the elliptic
equation
\begin{equation}\label{eq conduttivita'}
\mbox{div}(\sigma\nabla{u})=0\qquad\mbox{in}\quad\Omega,
\end{equation}
where the symmetric, positive definite, matrix $\sigma=\sigma(x)$,
$x\in\Omega$ represents the (possibly anisotropic) electric
conductivity. The inverse conductivity problem consists of finding
$\sigma$ when the so called Dirichlet-to-Neumann (D-N) map
$$
\Lambda_{\sigma}:u\vert_{\partial\Omega}\in{H}^{\frac{1}{2}}(\partial\Omega)
\longrightarrow{\sigma}\nabla{u}\cdot\nu\vert_{\partial\Omega}\in{H}^{-\frac{1}{2}}(\partial\Omega)
$$
is given for any $u\in{H}^{1}(\Omega)$ solution to (1.1). Here,
 $\nu$ denotes the unit outer normal to $\partial\Omega$. If
measurements can be taken only on one portion $\Gamma$ of
$\partial\Omega$, then the relevant map is called the local
Dirichlet-to-Neumann map. Let $\Gamma$ be a non-empty open portion
of $\partial\Omega$ and let us introduce the subspace of
$H^{\frac{1}{2}}(\partial\Omega)$
\begin{equation}\label{Hco}
H^{\frac{1}{2}}_{co}(\Gamma)=\big\{f\in
H^{\frac{1}{2}}(\partial\Omega) \:\vert\:\textnormal{supp}
\:f\subset\Gamma\big\}.
\end{equation}
The local Dirichlet-to-Neumann map is given, in a weak
formulation, as the operator $\Lambda_{\sigma}^{\Gamma}$ such that
\begin{equation}
\langle\Lambda_{\sigma}^{\Gamma}\;u,\;\phi\rangle=\int_{\Omega}\sigma\;\nabla
u\cdot\nabla\phi,
\end{equation}
for any $u,\;\phi\in{H}^{1}(\Omega)$,
$u\vert_{\partial\Omega},\;\phi\vert_{\partial\Omega}\in{H}_{co}^{\frac{1}{2}}(\Gamma)$
and $u$ is a weak solution to (1.1).

The problem of recovering the conductivity of a body by taking
measurements of voltage and current on its surface has came to be
known as Electrical Impedance Tomography (EIT). Different
materials display different electrical properties, so that a map
of the conductivity $\sigma(x)$, $x\in\Omega$ ($\Omega$ domain in
$\mathbb R^{n}$) can be used to investigate internal properties of
$\Omega$. EIT has many important applications in fields such as
geophysics, medicine and non--destructive testing of materials.
The first mathematical formulation of the inverse conductivity
problem
is due to A. P. Calder\`{o}n \cite{C}, 
 where he addressed
the problem of whether it is possible to determine the (isotropic)
conductivity by the D-N map. The case when measurements can be
taken all over the boundary has been studied extensively in the
past and fundamental papers like  \cite{KV1}, \cite{KV2},
 \cite{SU} and \cite{A} show that the isotropic case can be considered
solved. On the other hand the anisotropic case is still open and
different lines of research have been pursued. One direction has
been to find the conductivity up to a diffeomorphism which keeps
the boundary fixed (see \cite{LU}, \cite{S}, \cite{N} and
\cite{LaU}).

Another direction has been the one to assume that the anisotropic
conductivity is \textit{a priori} known to depend on a restricted
number of spatially--dependent parameters
(see \cite{A}, \cite{AG}, \cite{GL} and \cite{L}).\\
The problem of recovering the conductivity $\sigma$ by the
knowledge of the local Dirichlet-to-Neumann map
$\Lambda_{\sigma}^{\Gamma}$ has been treated more recently. Lassas
and Uhlmann in \cite{LaU} recovered a connected compact
real-analytic Riemannian manifold $(M,\:g)$ with boundary by
making use of the Green's function of the Laplace-Beltrami
operator $\Delta_{g}$. See also \cite{LaUT}.

The procedure of reconstructing the conductivity at the boundary by local
measurements has been studied first by Brown \cite{B}, where the
author gives a formula for reconstructing the isotropic
conductivity pointwise at the boundary of a Lipschitz domain
$\Omega$ without any \textit{a priori} smoothness assumption of
the conductivity. Nakamura and Tanuma \cite{NaT1} give a formula
for the pointwise reconstruction of a conductivity continuous at
one point $x^{0}$ of the boundary from the local D-N map when the
boundary is $C^{1}$ near $x^{0}$. Under some additional regularity
hypothesis the authors give a reconstruction formula for the
normal derivatives of $\sigma$ on $\partial\Omega$ at
$x^{0}\in\partial\Omega$ up to a certain order. A direct method
for reconstructing the normal derivative of the conductivity from
the local D-N map is presented in \cite{NaT2}. The result in
\cite{NaT1} has been improved by Kang and Yun \cite{KY} to an
inductive reconstruction method by using only the value of
$\sigma$ at $x^{0}$. The authors derive here also H\"older
stability estimates for the inverse problem to identify Riemannian
metrics (up to isometry) on the boundary via the local D-N map. An
overview on reconstructing formulas of the conductivity and its
normal derivative can be found in \cite{NaT3}.

For related results of uniqueness in the interior in the case  of local boundary
data, we refer to Bukhgeim and Uhlmann \cite{BU}, Kenig,
Sj\"ostrand and  Uhlmann \cite{ksu} and Isakov \cite{I}, and, for
stability, Heck and Wang \cite{hw}. Results of stability for cases
of piecewise constant conductivities and local boundary maps have
also been obtained by Alessandrini and Vessella \cite{AV} and by
Di Cristo \cite{D}.

It should also be emphasized that, mainly for the applications of
medical imaging, and in particular for breast cancer detection by
EIT, rather than the local Dirichlet-to-Neumann map, one should
consider the so-called local Neumann-to-Dirichlet (N-D) map. That
is, the map associating to specified current densities supported
on a portion $\Gamma\subset\partial\Omega$ the corresponding
boundary voltages, also measured on the same portion $\Gamma$ of
$\Omega$.

In the present paper we study the inverse conductivity problem by
local  maps, concentrating on the issue of determining the
boundary values of the conductivity and of its derivatives. We
continue the line of investigation pursued in \cite{AG}, by
considering anisotropic unknown conductivities having the
structure $\sigma(x) = A(x,a(x))$, where $A(x,t)$ is a known,
matrix valued, function and $a(x)$  is an unknown scalar function.
The precise assumptions shall be illustrated in section
\ref{results}. We improve upon the results obtained in
\cite{AG} under the following aspects.\\
(i) The uniqueness and stability estimates are adapted to the
local D-N map (Theorems \ref{stabilita' al bordo}, \ref{stabilita'
Holder al bordo}, \ref{unicita' al bordo} and
Corollary \ref{unicita' all'interno}),\\
(ii) the stability estimate at the boundary is obtained in the
wider class of conductivities which are merely continuous in a
neighborhood of some point at the boundary (Theorem \ref{stabilita' puntuale}),\\
(iii) analogous results are obtained when the local D-N map is
replaced by the local N-D map (Theorem \ref{stabilita' puntuale
N}).

 The paper is organized as follows. The main results are contained
in section \ref{results} (subsections \ref{D-to-N}, \ref{N-to-D}
for the local D-N, N-D maps respectively), while section \ref{sol
G zero fuori da Gamma} is devoted to the construction of singular
solutions of equation \eqref{eq conduttivita'} having the same
type of singularity as those in [A] but having support compactly
embedded on a non--empty open subset of the boundary (see Theorem
\ref{teore soluzioni singolari 2}). Proofs of the main results are
given in section \ref{proofs} (subsections \ref{D-to-N proofs},
\ref{N-to-D proofs}
for the local D-N, N-D maps respectively).\\
\section*{\normalsize{Acknowledgments}}
The authors gratefully acknowledge a fruitful conversation with D.
Isaacson who first suggested to consider the local
Neumann-to-Dirichlet (N-D) map for medical imaging applications.
This study was initiated when R.G. held a research contract at the
Dipartimento di Matematica ed Informatica at the University of
Trieste. R.G. wishes to acknowledge also the support of Science
Foundation Ireland (Grant 03/IN3/I401). The research of G.A. was
supported in part by MIUR (Grant  2006014115).
%
%
%
%
\section{\normalsize{Main results.}}\label{results} Let $\Omega$ be a domain in $\mathbb{R}^{n}$ $(n\geq 3)$,
with Lipschitz boundary $\partial\Omega$.  We recall, for sake of
completeness, the definition of Lipschitz regularity of the
boundary. We stick to notation already used in \cite{AG}.
\begin{defin}
Given positive numbers $L,\:r,\:h$ satisfying $h\geq{L}r$, we say
that a bounded domain $\Omega\in\mathbb{R}^{n}$ has Lipschitz
boundary if, for every $x^{0}\in\partial\Omega$, there exists a
rigid transformation of coordinates which maps $x^{0}$ into the
origin, such that, setting $x\:=\:(x^{\:\prime},\:x_{n})$,
$x^{\:\prime}\in\mathbb{R}^{n-1}$, $x_{n}\in\mathbb{R}$, we have
\begin{eqnarray*}
& &\Omega\cap\{\:x\:=\:(x^{\:\prime},\:x_{n})\:\vert
\quad\vert{x}^{\:\prime}\vert<\:r,\:\vert{x}_{n}\vert<\:h\}=\nonumber\\
& &=\:\{\:x\:=\:(x^{\:\prime},\:x_{n})\:\vert
\quad\vert{x}^{\:\prime}\vert<\:r,\:\vert{x}_{n}\vert<\:h,\:x_{n}\geq{f}(x^{\:\prime})\:\},
\end{eqnarray*}
\noindent{where} $f\:=\:f(x^{\:\prime})$ is a Lipschitz function
defined for $\vert{x}^{\:\prime}\vert<\:r$, which satisfies
\begin{eqnarray*}
f(0)&=&0\nonumber\\
\vert{f}(x^{\:\prime})-f(y^{\:\prime})\vert&\leq&L\:\vert{x}^{\:\prime}-y^{\:\prime}\vert,
\end{eqnarray*}
\noindent{for} every
$x^{\:\prime},\:y^{\:\prime}\in\mathbb{R}^{n-1}$, with
$\vert{x}^{\:\prime}\vert,\:\vert{\:y}^{\:\prime}\vert<\:r$.
\end{defin}
Let us now recall the class $\mathcal{H}$ of functions $A(x,\:t)$
introduced in \cite{AG}, which will be considered as admissible
conductivities.
\begin{defin}
Given $p>n$ , the positive constants
$\lambda,\:\mathcal{E},\:\mathcal{F}>0$, and denoting by $Sym_{n}$
the class of $n\times{n}$ real valued symmetric matrices, we say
that $A(\cdot,\:\cdot)\in\mathcal{H}$
if the following conditions hold\\
\begin{equation}
A\in{W}^{1,\:p}(\Omega\times{[\lambda^{-1},\:\lambda]\:,Sym_{n}}),
\end{equation}
\begin{equation}
D_{\:t}A\in{W}^{1,\:p}(\Omega\times{[\lambda^{-1},\:\lambda]}\:,\:Sym_{n}),
\end{equation}
\begin{eqnarray}
& &\textnormal{supess}_{\:t\in{\:[\lambda^{-1},\lambda]}}
\bigg(\parallel{A}(\cdot,t)\parallel_{L^{p}(\Omega)}+
\parallel{D}_{x}A(\cdot,t)\parallel_{L^{p}(\Omega)}\nonumber\\
& &+\parallel{D}_{t}A(\cdot,t)\parallel_{L^{p}(\Omega)}+
\parallel{D}_{t}D_{x}A(\cdot,t)\parallel_{L^{p}(\Omega)}\bigg)\leq\mathcal{E},
\end{eqnarray}
\begin{eqnarray}\label{ellitticita'}
\lambda^{-1}\vert\xi\vert^{2}\leq{A}(x,t)\xi\cdot\xi\leq\lambda\vert\xi\vert^{2},
& &for\:almost\:every\:x\in\Omega,\nonumber\\
&
&for\:every\:t\in{[\lambda^{-1},\lambda]},\:\xi\in\mathbb{R}^{n}.
\end{eqnarray}
\begin{eqnarray}\label{monotonia}
D_{t}A(x,t)\:\xi\cdot\xi\geq\mathcal{F}\vert\xi\vert^{2},
& & for\:almost\:every\:x\in\Omega\:,\nonumber\\
 & &
for\:every\:t\in{[\lambda^{-1},\lambda]}\:,\:
\xi\in\mathbb{R}^{n}.
\end{eqnarray}
\end{defin}
\noindent We observe that \eqref{ellitticita'} is a condition of
uniform ellipticity, whereas \eqref{monotonia} is a condition of
monotonicity with respect to the last variable $t$.
\begin{defin}\label{def rho}
For every $\rho$, $0<\rho<r$ we shall denote
\begin{eqnarray}
& &\Gamma_{\rho}=\big\{x\in\Gamma\vert\:\textnormal{dist}(x,\:\partial\Gamma)>\rho\big\},\\
& &U_{\rho}=\big\{x\in\mathbb R^{n}\vert\:\textnormal{dist}(x,\:
\Gamma_{\rho})<\frac{\rho}{4}\big\},\\
& &U_{\rho}^{i}=U_{\rho}\cap\Omega.
\end{eqnarray}
\end{defin}
\noindent Here it is understood that for the empty set
$\emptyset$, we have $\textnormal{dist}(x,\:\emptyset)=+\infty$.
It is evident that, $\Gamma$ being open and non--empty, there
exists $\rho_{0}$, $0<\rho_{0}\leq r$ such that
$\Gamma_{\rho_{0}}$ is also non empty. From now on we shall only
consider values of $\rho$ below $\rho_{0}$.
\begin{osse}
We emphasize that $\rho_{0}>0$ is a number which depends on the
choice of $\Gamma\subset\partial\Omega$. It should be evident that
if we choose $\Gamma$ narrower and narrower, then $\rho_{0}$ tends
to $0$ and one should expect a deterioration in the stability
estimates.
\end{osse}
\subsection{\normalsize{The Dirichlet-to-Neumann map.}}\label{D-to-N}
We start by rigorously defining the local D-N map. We consider a
given conductivity $\sigma\in L^{\infty}(\Omega\:,Sym_{n})$
satisfying the ellipticity condition
\begin{eqnarray}\label{ellitticita'sigma}
\lambda^{-1}\vert\xi\vert^{2}\leq{\sigma}(x)\xi\cdot\xi\leq\lambda\vert\xi\vert^{2},
& &for\:almost\:every\:x\in\Omega,\nonumber\\
& &for\:every\:\xi\in\mathbb{R}^{n},
\end{eqnarray}
and we fix an open, non--empty subset $\Gamma$ of
$\partial\Omega$. We denote by $<\cdot,\:\cdot>$ the
$L^{2}(\partial\Omega)$-pairing between
$H^{\frac{1}{2}}(\partial\Omega)$ and its dual
$H^{-\frac{1}{2}}(\partial\Omega)$.
\begin{defin}
The local Dirichlet-to-Neumann map associated to $\sigma$ and
$\Gamma$ is the operator
\begin{equation}\label{mappaDN}
\Lambda_{\sigma}^{\Gamma}:H^{\frac{1}{2}}_{co}(\Gamma)\longrightarrow
({H}^{\frac{1}{2}}_{co}(\Gamma))^{\ast}
\end{equation}
 defined by
\begin{equation}\label{def DN locale}
<\Lambda_{\sigma}^{\Gamma}\:g,\:\eta>\:=\:\int_{\:\Omega}
\sigma(x) \nabla{u}(x)\cdot\nabla\phi(x)\:dx,
\end{equation}
for any $g$, $\eta\in H^{\frac{1}{2}}_{co}(\Gamma)$, where
$u\in{H}^{1}(\Omega)$ is the weak solution to
\begin{displaymath}
\left\{ \begin{array}{ll} \textnormal{div}(\sigma(x)\nabla
u(x))=0, &
\textrm{$\textnormal{in}\quad\Omega$},\\
u=g, & \textrm{$\textnormal{on}\quad{\partial\Omega},$}
\end{array} \right.
\end{displaymath}
and $\phi\in H^{1}(\Omega)$ is any function such that
$\phi\vert_{\partial\Omega}=\eta$ in the trace sense.
\end{defin}
Note that, by \eqref{def DN locale}, it is easily verified that
$\Lambda^{\Gamma}_{\sigma}$ is selfadjoint. We shall denote by
$\parallel\cdot\parallel_{*}$ the norm on the Banach space of
bounded linear operators between $H^{\frac{1}{2}}_{co}(\Gamma)$
and $(H^{\frac{1}{2}}_{co}(\Gamma))^{*}$.
%
%

We can now state a first stability result for the boundary values of the conductivity, assuming that the unknown \emph{anisotropic} conductivity has the structure $\sigma(x) =A(x,\:a(x))$ with $A\in\mathcal H$ \emph{known} and $a=a(x)$ \emph{unknown scalar} function.
\begin{teore}\label{stabilita' al bordo}(\textnormal{Lipschitz stability of boundary
values}).
Given $p>n$, let $\Omega$ be a bounded domain with Lipschitz
boundary with constants L, $r$, h. Let $\Gamma$ be the subset of
$\partial\Omega$ introduced above and $\rho_{0}=\rho_{0}(\Gamma)$
the positive number introduced in Definition 2.3. If a, b are two
real-valued functions satisfying
\begin{equation}\label{limitazioni per a,b}
\lambda^{-1}\leq{a}(x),
b(x)\leq\lambda,\qquad{for}\:every\quad{x}\in\Omega,
\end{equation}
\begin{equation}\label{holder a,b}
\parallel{a}\parallel_{\:W^{1,p}(\Omega)}\:,\:\parallel{b}\parallel_{\:W^{1,p}(\Omega)}\:\leq{E},
\end{equation}
\noindent for some positive constant $E>0$ and $A\in\mathcal H$,
then for any $\rho$, $0<\rho\leq\rho_{0}$
\begin{equation}\label{stabilita' anisotropa}
\parallel{A}(x,\:a(x))-A(x,\:b(x))\parallel_{L^{\infty}\:(\:\bar\Gamma_{\rho})}
\leq{C}\parallel\Lambda_{A(x,\:a)}^{\Gamma}-\Lambda_{A(x,\:b)}^{\Gamma}\parallel_{*}.
\end{equation}
\noindent Here $C>0$ is a constant depending on n, p, L, r, h,
$diam(\Omega)$, $\rho$, $\rho_{0}$, $\lambda$, $\mathcal{E}$,
$\mathcal{F}$, but not on E.
\end{teore}
%
%

The next Theorem improves upon the previous one, in that the regularity assumption \eqref{holder a,b} is relaxed to mere continuity.
\begin{teore}\label{stabilita' puntuale}(\textnormal{Pointwise stability at the
boundary}). Given $p>n$, let $\Omega$, $\Gamma$ and $\rho_{0}$ be
as in \it Theorem \ref{stabilita' al bordo}. Suppose a, b are two
real valued functions satisfying \eqref{limitazioni per a,b} and
furthermore are continuous on $U_{\rho}^{i}$, for some $\rho$,
$0<\rho\leq\rho_{0}$. Let $A\in\mathcal{H}$, then for any
$x\in\Gamma_{\rho}$
\begin{equation}
\vert A(x,a(x))-A(x,b(x))\vert\leq C
\parallel\Lambda_{A(x,a)}^{\Gamma}-\Lambda_{A(x,b)}^{\Gamma}\parallel_{*},
\end{equation}
\noindent where $C>0$ is a constant which depends on n, p, L, r,
h, $diam(\Omega)$, $\rho_{0}$, $\rho$, $\lambda$, $\mathcal{E}$,
$\mathcal{F}$ only.
\end{teore}
%
%

Here we state our stability results for boundary values of the derivatives of the conductivity.
\begin{teore}\label{stabilita' Holder al bordo}(\textnormal{H\"older stability of derivatives at the
boundary}). Given $p$, $\Omega$, $\Gamma$ and $\rho_{0}$ as in
Theorem 2.1, let a, b satisfy \textnormal{(2.11)},
\textnormal{(2.12)} and $A\in\mathcal{H}$. Suppose furthermore
that for some $\rho$, $0<\rho\leq\rho_{0}$, some positive integer
k and some $\alpha$, $0<\alpha<1$ we have
\begin{equation}
A\in{C}^{\:k,\:\alpha}(\:\bar{U}_{\rho}\times{[\lambda^{-1},\lambda]\:}),
\end{equation}
\begin{equation}
\parallel{A}\parallel_{\:C^{\:k,\:\alpha}(\:\bar{U}_{\rho}\times{[\lambda^{-1},\:\lambda]})}\:\leq{E}_{\:k}.
\end{equation}
\begin{equation}
\parallel{a-b}\parallel_{\:C^{\:k,\:\alpha}(\:\bar{U}_{\rho})}\:\leq{E}_{\:k}.
\end{equation}
\noindent Then
\begin{eqnarray}
&
&\parallel{D}^{\:k}(A(x,a(x))-A(x,b(x)))\parallel_{L^{\:\infty}\:
(\bar\Gamma_{\rho})}\nonumber\\
& &\leq{C}\parallel\Lambda_{A(x,\:a)}^{\Gamma}-
\Lambda_{A(x,\:b)}^{\Gamma}\parallel_{*}^{\:\delta_{\:k}\:\alpha},
\end{eqnarray}
\noindent where
\begin{equation}
\delta_{k}\:=\:\prod_{j\:=\:0}^{k}\frac{\alpha}{\alpha+j}.
\end{equation}
\noindent Here $C>0$ is a constant which depends only on n, p, L,
r, h, $diam(\Omega)$, $\rho_{0}$, $\rho$, $\lambda$, E, $\alpha$,
k, and $E_{k}$.
\end{teore}
%
%

Under a slightly weaker assumption, we can also obtain the
following uniqueness result.
\begin{teore}\label{unicita' al bordo}(\textnormal{Uniqueness at the boundary}).
Let $p$, $\Omega$, $\Gamma$, $\rho_{0}$, $a$, $b$, $A$ as in
Theorem 2.3. Suppose that for some $\rho$, $0<\rho\leq\rho_{0}$
and some positive integer k we have
\begin{equation}
a-b\in{C}^{\:k}(\bar{U}_{\rho}).
\end{equation}
\noindent If
\begin{displaymath}
\Lambda_{A\:(x,\:a(x))}^{\Gamma}\:=\:\Lambda_{A\:(x,\:b(x))}^{\Gamma},
\end{displaymath}
\noindent then
\begin{equation}
D^{j}(a-b)\:=\:0\qquad{on}\quad\Gamma_{\rho},
\qquad{for}\:all\quad{j}\leq{k}.
\end{equation}
\noindent If in addition we have
\begin{equation}
A\in{C}^{\:k}\Big(
\bar{U}_{\rho}\times{[\lambda^{-1},\lambda]}\Big),
\end{equation}
\noindent{then}
\begin{equation}
D^{j}\bigg(A(x\:,\:a(x))\:=\:A(x\:,\:b(x))\bigg)=0
\qquad{on}\quad\Gamma_{\rho}, \quad{for}\:all\quad{j}\leq{k}.
\end{equation}
\end{teore}
%
%

What follows is a well--known consequence of the previous Theorem, see \cite{KV2} and \cite{A} for related arguments.
\begin{corol}\label{unicita' all'interno}(\textnormal{Uniqueness in the interior}).
Let a, b satisfy (2.11), (2.12) with $p\:=\:\infty$. Let
$A\in\mathcal{H}$ and in addition $A\in{W}^{1,\:\infty}\big(
\Omega\times{[\lambda^{-1},\lambda]},\:Sym_{n}\big) $. Suppose
that $\Omega$ can be partitioned into a finite number of Lipschitz
domains,
$\{A_{j}\}_{j\:=\:1,\ldots,N}$, such that $a-b$ is analytic on each $\bar{A}_{j}$.\\
\indent{If}
\begin{displaymath}
\Lambda_{A(x,\:a)}^{\Gamma}\:=\:\Lambda_{A(x,\:b)}^{\Gamma},
\end{displaymath}
\noindent then we have
\begin{equation}
A(x,\:a(x))\:=\:A(x,\:b(x))\qquad{in}\quad\Omega.
\end{equation}
\end{corol}
\subsection{\normalsize{The Neumann-to-Dirichlet map.}}\label{N-to-D}
Let us introduce the following function spaces
\begin{equation*}
_{0}H^{\frac{1}{2}}(\partial \Omega)=\left\{\phi\in
H^{\frac{1}{2}}(\partial \Omega)\vert\:\int_{\partial\Omega}\phi\:
=0\right\},
\end{equation*}
\begin{equation*}
_{0}H^{-\frac{1}{2}}(\partial \Omega)=\left\{\psi\in
H^{-\frac{1}{2}}(\partial \Omega)\vert\:\langle\psi,\:1\rangle=0
\right\}.
\end{equation*}
Observe that if we consider the (global) D-N map
$\Lambda_{\sigma}$, that is the map $\Lambda_{\sigma}^{\Gamma}$ ,introduced in \eqref{mappaDN},
in the special case when $\Gamma =
\partial \Omega$, we have that it maps onto $_{0}H^{-\frac{1}{2}}(\partial \Omega)$, and, when restricted to $_{0}H^{\frac{1}{2}}(\partial
\Omega)$, it is injective with bounded inverse. Then we can define
the global Neumann-to-Dirichlet map as follows.
\begin{defin}
The  Neumann-to-Dirichlet map associated to $\sigma$, $N_{\sigma}
:\  _{0}H^{-\frac{1}{2}}(\partial \Omega)\longrightarrow \
_{0}H^{\frac{1}{2}}(\partial \Omega)$ is given by
\begin{equation}\label{Nsigma}
N_{\sigma} = \left(\Lambda_{\sigma}|_{_{0}H^{\frac{1}{2}}(\partial
\Omega)} \right)^{-1} .\end{equation}
\end{defin}
Note that $N_{\sigma}$ can also be characterized as the selfadjoint
operator satisfying
\begin{equation}\label{ND globale}
<\psi,\:N_{\sigma}\psi>\:=\:\int_{\:\Omega} \sigma(x)
\nabla{u}(x)\cdot\nabla{u}(x)\:dx,
\end{equation}
for every $\psi\in\: _{0}H^{-\frac{1}{2}}(\partial \Omega)$, where
$u\in{H}^{1}(\Omega)$ is the weak solution to the Neumann problem
\begin{equation}\label{N bvp}
\left\{ \begin{array}{lll} \textnormal{div}(\sigma\nabla u)=0, &
\textrm{$\textnormal{in}\quad\Omega$},\\
\sigma\nabla u\cdot\nu\vert_{\partial\Omega}=\psi, &
\textrm{$\textnormal{on}\quad{\partial\Omega}$},\\
\int_{\partial\Omega}u\: =0.
\end{array} \right.
\end{equation}
We are now in position to introduce the local version of the N-D
map. Let  $\Gamma$ be an open portion of $\partial\Omega$ and let
$\Delta=\partial\Omega\setminus\bar\Gamma$.
We denote by
$H^{\frac{1}{2}}_{00}(\Delta)$  the closure in
$H^{\frac{1}{2}}(\partial\Omega)$ of the space
$H^{\frac{1}{2}}_{co}(\Delta)$ previously defined in \eqref{Hco}
and we introduce
\begin{equation}
_{0}H^{-\frac{1}{2}}(\Gamma)=\left\{\psi\in \:
_{0}H^{-\frac{1}{2}}(\partial\Omega)\vert\:\langle\psi,\:f\rangle=0,\quad\textnormal{for\:any}\:f\in
H^{\frac{1}{2}}_{00}(\Delta)\right\},
\end{equation}
that is the space of distributions $\psi \in
H^{-\frac{1}{2}}(\partial\Omega)$ which are supported in
$\bar\Gamma$ and have zero average on $\partial\Omega$. The local
N-D map is then defined as follows.
\begin{defin}
The local Neumann-to-Dirichlet map associated to $\sigma$,
$\Gamma$ is the operator $N_{\sigma}^{\Gamma}:\:
_{0}H^{-\frac{1}{2}}(\Gamma)\longrightarrow
\big(_{0}H^{-\frac{1}{2}}(\Gamma)\big)^{\ast}\subset{_{0}H}^{\frac{1}{2}}(\partial\Omega)$
given  by
\begin{equation}
\langle N_{\sigma}^{\Gamma}\;i,\;j\rangle=\langle
N_{\sigma}\;i,\;j\rangle,
\end{equation}
for every $i, j\in\: _{0}H^{-\frac{1}{2}}(\Gamma)$.
\end{defin}
When the local D-N map is replaced by the above defined local N-D
map, completely analogous results to Theorems \ref{stabilita' al
bordo}-\ref{unicita' al bordo} and Corollary \ref{unicita'
all'interno} could be obtained. For the sake of simplicity we
state the appropriate version of Theorem \ref{stabilita' puntuale}
only. See also Remark \ref{remark finale} for further details. In
what follows, we shall denote by
$\parallel\cdot\parallel_{\ast\:\ast}$ the norm on the Banach
space of bounded linear operators between
$_{0}H^{-\frac{1}{2}}(\Gamma)$ and
$\left(_{0}H^{-\frac{1}{2}}(\Gamma)\right)^{\ast}$.
\begin{teore}\label{stabilita' puntuale N}
Given $p>n$, let $\Omega$, $\Gamma$ and $\rho_{0}$ be as in \it
Theorem \ref{stabilita' al bordo}. Suppose a, b are two real
valued functions satisfying \eqref{limitazioni per a,b},
continuous on $U_{\rho}^{i}$, for some $\rho$,
$0<\rho\leq\rho_{0}$. Let $A\in\mathcal{H}$, then for any
$x\in\Gamma_{\rho}$
\begin{equation}
\vert A(x,a(x))-A(x,b(x))\vert\leq C
\parallel N_{A(x,\:a)}^{\Gamma}-N_{A(x,\:b)}^{\Gamma}\parallel_{\ast\:\ast},
\end{equation}
\noindent where $C>0$ is a constant which depends on n, p, L, r,
h, $diam(\Omega)$, $\rho_{0}$, $\rho$, $\lambda$, $\mathcal{E}$,
$\mathcal{F}$ only.
\end{teore}

%
%
\section{\normalsize{Singular solutions vanishing on
$\partial\Omega\setminus\Gamma$.}}\label{sol G zero fuori da
Gamma}
This section is devoted to the construction of particular
solutions of equation \eqref{eq conduttivita'}, having the same
type of singularity of those constructed in \cite {A} but
vanishing on the portion of the boundary
$\partial\Omega\setminus\Gamma$. We consider the elliptic operator
\begin{equation}\label{operatore L}
L\:=\:\frac{\partial}{\partial{x}_{i}}\bigg(
\sigma_{ij}\frac{\partial}{\partial{x}_{j}}\bigg) ,
\qquad\textnormal{in}\quad{B}_{R}\:=\:\big\{x\in\mathbb{R}^{n}\vert\:\vert{x}\vert<R\big\},
\end{equation}
\noindent where the coefficient matrix $(\sigma_{ij}(x))$ is
symmetric and satisfies
\begin{equation}\label{ellitticita' operatore L}
\lambda^{-1}\vert\:\xi\vert^{2}\leq\sigma_{ij}(x)\xi_{i}\:\xi_{j}\leq\lambda\vert\:\xi\vert^{2},
\qquad{\textnormal{for\:every}}\quad{x}\in{B}_{R},\:\xi\in\mathbb{R}^{n},
\end{equation}
\noindent and also
\begin{equation}\label{norma W1p}
\parallel\sigma_{ij}\parallel_{\:W^{1,\:p}(B_{R})}\leq{E},\qquad{i},\:j\:=\:1,\ldots,n,
\end{equation}
\noindent here $p>n$ and $\lambda,\:E$ are positive constants. We
recall the following theorem from \cite{A}.
\begin{teore}\label{teore soluzioni singolari}(\textnormal{Singular Solutions}).
Let L satisfy \eqref{operatore L}-\eqref{norma W1p}. For every
spherical harmonic $S_{m}$ of degree $m\:=\:0,\:1,\:2,\:...,$
there exists $u\in{W}^{2,\:p}_{loc}(B_{R}\setminus\{0\})$ such
that
\begin{equation}
Lu\:=\:0\qquad{in}\quad{B}_{R}\setminus\{\:0\},
\end{equation}
\noindent and furthermore,
\begin{eqnarray}\label{soluzioni singolari 1}
u(x)&=&\log\:\vert Jx\vert\:S_{\:0}\Bigg(
\frac{Jx}{\vert{\:Jx}\:\vert}\Bigg) +w(x),
\:when\:n\:=\:2\:and\:m=0,\nonumber\\
u(x)&=&\vert Jx\vert^{\:2-n-m}\:S_{m}\Bigg(
\frac{Jx}{\vert{\:Jx}\:\vert}\Bigg) +w(x), \:otherwise,
\end{eqnarray}
\noindent where J is the positive definite symmetric matrix such
that $J\:=\:\sqrt{(\sigma_{ij}(0))^ {-1}}$ and w satisfies
\begin{equation}\label{w stima lipschitz}
\vert{\:w}(x)\vert+\vert{\:x}\:\vert\:\vert{D}w(x)\vert\leq{C}\:\vert{\:x}\:\vert^
{\:2-n-m+\alpha},\quad{in}\quad{B}_{R}\setminus\{\:0\:\},
\end{equation}
\begin{equation}\label{w stima int}
\bigg( \int_{s<\vert{x}\vert<2s} \vert{D}^{2}w\vert^{p}\bigg) ^
{\frac{1}{p}}\leq{C}\:s^{-n-m+\alpha+\frac{n}{p}},\quad{for}\:every\quad{s},\:0<s<R/2.
\end{equation}
\noindent Here $\alpha$ is any number such that
$0<\alpha<1-\frac{n}{p}$, and C is a constant depending only on
$\alpha,\:n,\:p,\:R,\:\lambda,\:and\:E$.
\end{teore}
\textit{Proof}. See \cite[Theorem 1.1]{A}. $\quad\blacksquare$\\
We shall also need the following.
\begin{lemma}
Let the hypotheses of Theorem \ref{teore soluzioni singolari} be
satisfied. For every $m\:=\:0,\:1,\:2,\ldots$ there exists a
spherical harmonic $S_{m}$ of degree $m$ such that the solution
$u$ given by Theorem \ref{teore soluzioni singolari} also
satisfies
\begin{equation}
\vert{D}u(x)\vert>\vert{x}\vert^{\:1-(n+m)},\qquad
for\:every\:x,\quad 0<\vert{x}\vert<r_{0},
\end{equation}
\noindent where $r_{0}>0$ depends only on $\lambda$, E, p, m and
R.
\end{lemma}
\textit{Proof}. The proof of this lemma can be obtained along the
same lines as of \cite[Lemma 3.1]{A} and \cite[Section 3]{AG}.
$\quad\blacksquare$

%
%
Let us construct now solutions $u$ of \eqref{eq conduttivita'}
having a singularity of the same type of the above theorem in a
point outside $\Omega$ and satisfying
\begin{equation*}
u\vert_{\partial\Omega}\in{H^{\frac{1}{2}}_{co}}(\Gamma),
\end{equation*}
\noindent in the sense of traces. To this purpose we shall make
use of an augmented domain $\Omega_{\rho}$. In fact, for any
$\rho$, $0<\rho\leq\rho_{0}$, one can always construct a domain
$\Omega_{\rho}$ with Lipschitz constants depending only on $\rho$,
$r$, $L$, $h$ such that
\begin{equation}
\Omega\subset\Omega_{\rho},\quad\partial\Omega\cap\Omega_{\rho}\subset\subset\Gamma
\end{equation}
and
\begin{equation}
\textnormal{dist}(x,\:\partial\Omega_{\rho})\geq\frac{\rho}{2},\quad\textnormal{for}\:\textnormal{every}
\quad x\in U_{\rho}.
\end{equation}
\noindent If $L$ is an operator of type \eqref{operatore L} on
$\Omega$, satisfying \eqref{ellitticita' operatore L},
\eqref{norma W1p} on $\Omega$, then for any $\rho>0$, one can
always extend the operator $L$ to $\Omega_{\rho}$ in such a way so
that $L$ still satisfies \eqref{ellitticita' operatore L},
\eqref{norma W1p} on the enlarged domain $\Omega_{\rho}$.
%
%
As the boundary $\partial\Omega$ is Lipschitz the unit normal
vector field to the boundary may not be defined pointwise so we
shall introduce a unitary vector field $\tilde\nu$ locally defined
near $\partial\Omega$ such that: (i) $\tilde\nu$ is $C^{\:\infty}$
smooth, (ii) $\tilde\nu$ is non-tangential to $\partial\Omega$
(see \cite[Section 3]{AG}, for the construction procedure of the
latter). The point $z_{\tau}=x^{0}+\tau\tilde\nu$, where
$x^{0}\in\partial\Omega$, satisfies
\begin{equation}\label{distanza z}
C\tau\leq d(z_{\tau},\:\partial\Omega)\leq\tau,
\end{equation}
\noindent for any $\tau$, $0<\tau\leq\tau^{0}$. Here $C$ and
$\tau^{0}$ are positive constants depending only on $L,\:r,\:h$
\cite[Lemma 3.3]{AG}.\\
We distinguish the cases when $m=0$ or $m>0$. For the case $m=0$
we recall the following asymptotic estimate which only requires
the H\"older continuity of the coefficients.
\begin{teore}\label{teore funzione green}
Let $\Omega$ and $\Gamma$ be as in Theorem \ref{stabilita' al
bordo}. For any $\tau$, $0<\tau\leq\tau_{0}$, set
$z_{\tau}=x^{0}+\tau\tilde\nu$, for some
$x^{0}\in\bar\Gamma_{\rho}$ and $\rho$, $0<\rho\leq\rho_{0}$. If
$L$ is the operator  of \eqref{operatore L}, with H\"older
continuous coefficients matrix
$\sigma=\{\sigma_{ij}\}_{i,\:j=1\dots n}$, with exponent
$0<\beta<1$, the Green's function $G_{\sigma}$ for the Dirichlet
boundary value problem
\begin{displaymath}
\left\{ \begin{array}{ll}
L\:G_{\sigma}(x,\:z_{\tau})=-\delta(x-z_{\tau}), & \textrm{$\textnormal{in}\quad\Omega_{\rho}$}\\
G_{\sigma}(\cdot,\:z_{\tau})=0, &
\textrm{$\textnormal{on}\quad{\partial\Omega_{\rho}}$}
\end{array} \right.
\end{displaymath}
has the form
\begin{equation}\label{espression greens function}
G_{\sigma}(x,\:z_{\tau})=
C_{n}\:\big(\det(\sigma(z_{\tau}))\big)^{-1/2}\Big(\sigma^{-1}(z_{\tau})(x-z_{\tau})\cdot(x-z_{\tau})\Big)^{\frac{2-n}{2}}+R(x,\:z_{\tau}),
\end{equation}
\noindent where $C_{n}$ is a suitable dimensional constant and the
remainder $R(x,\:z_{\tau})$ satisfies
\begin{equation}\label{stima su R1}
\vert R(x,\:z_{\tau})\vert+\vert x-z_{\tau}\vert\:\vert
\nabla_{x}R(x,\:z_{\tau})\vert \leq C \vert
x-z_{\tau}\vert^{2-n+\alpha},
\end{equation}
\noindent for every $x\in\Omega_{\rho}$, $\vert
x-z_{\tau}\vert\leq r_{0}$, where $C=C(E)$ is a positive constant
depending on $E$, $r_{0}$ is a positive number which depends only
on the geometry of $\Omega$ and $0<\alpha<\beta$.
\end{teore}
\textit{Proof of Theorem \ref{teore funzione green}}. We refer to
\cite[Chapter 1]{Mi} and \cite[(1.31)-(1.33)]{MT}. $\quad\blacksquare$\\
As a Corollary, we also have
\begin{corol}\label{corollario sol sing 1}
The Green's function $G_{\sigma}$ introduced in Theorem \ref{teore
funzione green} satisfies
\begin{equation}\label{stima H1}
\|\:G_{\sigma}(\cdot,\:z_{\tau})\|_{H^{1}(\Omega)}\leq
C\tau^{(2-n)/2},\qquad for\:any\quad0<\tau\leq\tau^{0},
\end{equation}
\noindent where $C>0$ is a constant which only depends on
$diam({\Omega})$, $\lambda$, L, r, h and $\tau^{0}$.
\end{corol}
\textit{Proof of Corollary \ref{corollario sol sing 1}}. A
straightforward consequence of the pointwise upper bound
\eqref{distanza z} and of the  Caccioppoli Inequality (see, for example,
\cite[Chapter 7 ]{Gi}) yields
\begin{equation}\label{Caccioppoli}
\|\:G_{\sigma}\|_{H^{1}(\Omega)}\leq\:\frac{K}{\tau}\:\|\:G_{\sigma}\|_{L^{2}(\Omega_{\rho}\setminus
B_{\:C\frac{\tau}{2}}(z_{\tau}))},
\end{equation}
\noindent where $K=K(\lambda,\:L,\:r,\:h,\:\tau^{0})$ is a
positive constant depending only on $\lambda$, $L$, $r$, $h$ and
$\tau^{0}$ and $C>0$ in \eqref{Caccioppoli} is the constant
introduced in \eqref{distanza z}.$\quad\blacksquare$\\
For the case $m>0$, we shall need stronger regularity assumptions
on the coefficients. In fact, under the $W^{1,\:p}$ bound
\eqref{norma W1p} we obtain
\begin{teore}\label{teore soluzioni singolari 2}
Let $\Omega$ and $\Gamma$ be as in Theorem \ref{stabilita' al
bordo} For any $\rho$, $0<\rho\leq\rho_{0}$, let z be an arbitrary
point in $U_{\rho}$. For every $m=0,\:1,\:2,\dots$ and for every
spherical harmonic $S_{m}\neq 0$ of degree m, there exists $u\in
H^{1}_{loc}(\overline\Omega_{\rho}\setminus\{z\}) \cap
W^{2,\:p}_{loc}(\Omega_{\rho}\setminus\{z\})$ such that
\begin{equation}\label{eq teore sol sing 2}
Lu=0\qquad{in}\quad\Omega_{\rho}\setminus\{z\},
\end{equation}
\begin{equation}\label{cond bordo teore sol sing 2}
u=0\quad on\quad\partial\Omega_{\rho},\quad in\:the\:trace\:sense
\end{equation}
\noindent and it has the form
\begin{equation}\label{soluzioni singolari 2}
u(x)=\vert J(x-z)\vert^{\:2-n-m}\:S_{m}\Bigg(
\frac{J(x-z)}{\vert{\:J(x-z)}\:\vert}\Bigg) +v(x),
\end{equation}
\noindent{where} J is the positive definite symmetric matrix such
that $J\:=\:\sqrt{(\sigma_{ij}(z))^ {-1}}$ and the remainder $v$
satisfies
\begin{equation}\label{v stima lipschitz}
\vert{\:v}(x)\vert+\vert{\:x-z}\:\vert\:\vert{D}v(x)\vert\leq{C}\:\vert{\:x-z}\:\vert^
{\:2-n-m+\alpha},\quad{in}\quad{B}_{\rho/4}(z)\setminus\{\:z\:\},
\end{equation}
\begin{equation}\label{v stima int}
\bigg( \int_{s<\vert{x}-z\vert<2s} \vert{D}^{2}v\vert^{p}\bigg) ^
{\frac{1}{p}}\leq{C}\:s^{-n-m+\alpha+\frac{n}{p}},\quad{for}\:every\quad{s},\:0<2s<\rho/4.
\end{equation}
\noindent Here $\alpha$ is any number such that
$0<\alpha<1-\frac{n}{p}$, and C is a constant depending only on
$\alpha,\:n,\:p,\:R,\:\lambda,\:\rho_{0},\:\rho\:and\:E$.
\end{teore}
\begin{osse}
Notice that, if $z\in U_{\rho}\setminus\Omega$ then $u\in
H^{1}(\Omega)$ and its trace satisfies $u\vert_{\partial\Omega}\in
H^{\frac{1}{2}}_{co}(\Gamma)$.
\end{osse}
\textit{Proof of Theorem \ref{teore soluzioni singolari 2}}. With
no loss of generality we can assume $z=0$. Consider a positive
number $R$ sufficiently large so that $B_{R}(0)$, the ball with
centre $0$ and radius $R$, is such that
$\Omega_{\rho}\subset B_{R/2}(0)\subset\subset B_{R}(0)$. We
consider the singular solution of Theorem \ref{teore soluzioni
singolari} on $B_{R}(0)$. Let us denote this solution by $u_{m}$.
Let $w_{0}$ be the solution to the problem
\begin{displaymath}
\left\{ \begin{array}{ll}
\mbox{div}(\sigma\nabla w_{0})=0, & \textrm{$\textnormal{in}\quad\Omega_{\rho}$}\\
w_{0}=-u_{m}, &
\textrm{$\textnormal{on}\quad{\partial\Omega_{\rho}}$.}
\end{array} \right.
\end{displaymath}
\noindent By recalling \eqref{soluzioni singolari 1} we get
\begin{equation}
\sup_{\partial\Omega_{\rho}} \big(\vert u_{m}\vert+\vert\nabla
u_{m}\vert\big)\leq C_{1},
\end{equation}
\noindent consequently
\begin{equation}
\parallel w_{0} \parallel_{H^{1}(\Omega_{\rho})} \leq C_{2},
\end{equation}
\noindent where $C_{1}$ is a positive constant which depends on
$\rho_{0}$, $\rho$, $n$ and $m$ only and $C_{2}>0$ depends only on
$\rho_{0}$, $\rho$, $n$, $m$, $R$, $L$, $r$ and $h$. If we set
\begin{equation}
u(x)=\vert Jx\vert^{\:2-n-m}\:S_{m}\Bigg(
\frac{Jx}{\vert{\:Jx}\:\vert}\Bigg)+v(x),
\qquad\textnormal{for}\:\textnormal{any}\quad x\in\Omega_{\rho}
\end{equation}
and $v=w+w_{0}$, where $w$ is the reminder appearing in
\eqref{soluzioni singolari 1}. Then $u$ can be written as
\begin{equation}
u(x)=u_{m}(x)+w_{0}(x),\qquad
\textnormal{for}\:\textnormal{any}\quad x\in\Omega_{\rho}
\end{equation}
and satisfies \eqref{eq teore sol sing 2}, \eqref{cond bordo teore
sol sing 2}, moreover, by a standard interior regularity estimate
\begin{equation}
\parallel w_{0} \parallel_{W^{2,\:p}(B_{\rho/4}(z))} \leq C,
\end{equation}
where $C>0$ depends on $\rho_{0}$, $\rho$, $m$ and $n$, $R$, $L$,
$r$ and $h$. Hence, recalling the bounds \eqref{w stima
lipschitz}, \eqref{w stima int} we obtain for $v=w+w_{0}$,
\eqref{v stima lipschitz}, \eqref{v stima int}.
$\quad\blacksquare$
%
%
\section{\normalsize{Proofs of the main theorems.}}\label{proofs}
\subsection{\normalsize{The D-N map.}}\label{D-to-N proofs}
The proofs of Theorem \ref{stabilita' Holder al bordo},
\ref{unicita' al bordo} and Corollary \ref{unicita' all'interno}
follow the same line of the corresponding results in \cite{AG} by
replacing the singular solutions used there by those introduced in
the previous Section 3 which vanish on
$\partial\Omega\setminus\Gamma$. For this reason, we shall give
the details of the proof of Theorem \ref{stabilita' al bordo}
only.

%
%
%
%
\textit{Proof of Theorem \ref{stabilita' al bordo}}.  Let
$x^{0}\!\!\in\!\bar\Gamma_{\rho}$ such that
$(a\!-\!b)(x^{0})=\parallel
\!\!a\!-\!b\!\!\parallel_{L^{\infty}(\Gamma_{\rho})}$ and set
$z_{\tau}=x^{0}+\tau\tilde\nu$, with
$0<\tau\leq\mbox{min}\big\{\tau_{0},\:\frac{\rho}{8}\big\}$. Let
$G_{a}$, $G_{b}$ be the Green's functions of Theorem \ref{teore
funzione green} in $\Omega_{\rho}$ for the operators
$\mbox{div}\left(A(\cdot,a(\cdot))\nabla\cdot\right)$ and
$\mbox{div}\left(A(\cdot,b(\cdot))\nabla\cdot\right)$
respectively, that is, for instance
\begin{displaymath}
\left\{ \begin{array}{ll}
\mbox{div}(A(x,a(x))\nabla G_{a}(x,\:z_{\tau}))=-\delta(x-z_{\tau}), & \textrm{$\textnormal{in}\quad\Omega_{\rho}$}\\
G_{a}(\cdot,\:z_{\tau})=0, &
\textrm{$\textnormal{on}\quad{\partial\Omega_{\rho}}$.}
\end{array} \right.
\end{displaymath}
\noindent and analogously for $G_{b}$. By possibly reducing $\tau$
and taking
$0<\tau\leq\mbox{min}\big\{\tau_{0},\:\frac{\rho}{8},\:\frac{r_{0}}{2}\big\}$,
we have that $B_{r_{0}}(z_{\tau})\cap\Omega$ is not empty and
moreover $B_{r_{0}}(z_{\tau})\cap\Omega\subset U_{\rho}$. By
recalling \eqref{def DN locale} and \cite[(b), p. 253]{A} we can
write
\begin{eqnarray}\label{Aless id con G}
 \int_{B_{r_{0}}(z_{\tau})\cap\Omega}
\!\!\!\!\!\!\!\!\!\!\!\!\!\!\!\!\big(A(x,a)-A(x,b)\big)\:\nabla
G_{a}\cdot\nabla G_{b} &\!\!\!+\!\!\!&\int_{\Omega\setminus
B_{r_{0}}(z_{\tau})}
\!\!\!\!\!\!\!\!\!\!\!\!\!\!\!\!\big(A(x,a)-A(x,b)\big)\:\nabla G_{a}\!\cdot\!\nabla G_{b}\nonumber\\
&\!\!\!=\!\!\!& \langle\:
(\Lambda_{A(x,a)}^{\Gamma}-\Lambda_{A(x,b)}^{\Gamma})\:G_{a},\:G_{b}\:\rangle.
\end{eqnarray}
\noindent Here and in the sequel, it is understood
$G_{a}=G_{a}(\cdot,\:z_{\tau})$ and analogously for $G_{b}$. By
combining \eqref{Aless id con G} with \eqref{espression greens
function} and \eqref{stima su R1}, we obtain
\begin{eqnarray}\label{stime Ja Jb}
& &\frac{(2-n)^{2}}{\big(\det A(a,\:z_{\tau})\big)^{1/2}\big(\det
A(b,\:z_{\tau})\big)^{1/2}}\nonumber\\
&\times&\int_{B_{r_{0}}(z_{\tau})\cap\Omega}
\frac{J_{b}^{2}\big(A(x,\:a)-A(x,\:b)\big)J_{a}^{2}\:(x-z_{\tau})\cdot(x-z_{\tau})}
{\big(J_{a}^{2}(x-z_{\tau})\!\cdot\!(x-z_{\tau})\big)^{n/2}\big(J_{b}^{2}(x-z_{\tau})\cdot(x-z_{\tau})\big)^{n/2}}\nonumber\\
&\leq& C(E,\:\mathcal{E})
\bigg\{\int_{B_{r_{0}}(z_{\tau})\cap\Omega} \vert
x-z_{\tau}\vert^{2(1-n)+\beta}\nonumber\\
&+&\int_{\Omega\setminus B_{r_{0}}(z_{\tau})}
\vert A(x,a)-A(x,b)\vert\vert x-z_{\tau}\vert^{2-2n}\bigg\}\nonumber\\
&+&
\parallel\Lambda_{A(x,\:a)}^{\Gamma}-\Lambda_{A(x,\:b)}^{\Gamma}\parallel_{*}\:
\parallel\! G_{a}\!\parallel_{H^{\frac{1}{2}}_{co}(\Gamma)}\:
\parallel\! G_{b}\!\parallel_{H^{\frac{1}{2}}_{co}(\Gamma)},
\end{eqnarray}
\noindent with
$J_{a}=\sqrt{\big(A(z_{\tau},a(z_{\tau}))\big)^{-1}}$,
$J_{b}=\sqrt{\big(A(z_{\tau},b(z_{\tau}))\big)^{-1}}$ and by the
H\"older continuity of $A(x,a(x))$, $A(x,b(x))$, \eqref{stime Ja
Jb} yields
\begin{eqnarray}\label{stima A(-1)}
& &C(n,\:\mathcal{E})\int_{B_{r_{0}}(z_{\tau})\cap\Omega}
\frac{\big(A(x^{0},\:a(x^{0}))^{-1}-A(x^{0},\:b(x^{0}))^{-1}\big)\:(x-z_{\tau})\cdot(x-z_{\tau})}
{\big(J_{a}^{2}(x-z_{\tau})\!\cdot\!(x-z_{\tau})\big)^{n/2}\big(J_{b}^{2}(x-z_{\tau})\cdot(x-z_{\tau})\big)^{n/2}}\nonumber\\
&\leq& C(E,\:\mathcal{E})
\bigg\{\int_{B_{r_{0}}(z_{\tau})\cap\Omega}
\vert x^{0}-z_{\tau}\vert^{\beta}\vert x-z_{\tau}\vert^{2(1-n)}\nonumber\\
\!\!\!&+&\!\!\!\!\int_{B_{r_{0}}(z_{\tau})\cap\Omega}
\vert x-x^{0}\vert^{\beta}\vert x-z_{\tau}\vert^{2(1-n)}\nonumber\\
\!\!\!&+&\!\!\!\! \int_{\Omega\setminus B_{r_{0}}(z_{\tau})}
\vert a-b\vert\vert x-z_{\tau}\vert^{2-2n}\bigg\}\nonumber\\
\!\!\!&+&\!\!\!\!
\parallel\Lambda_{A(x,\:a)}^{\Gamma}-\Lambda_{A(x,\:b)}^{\Gamma}\parallel_{*}
\parallel\! G_{a}\!\parallel_{H^{\frac{1}{2}}_{co}(\Gamma)}
\parallel\! G_{b}\!\parallel_{H^{\frac{1}{2}}_{co}(\Gamma)},
\end{eqnarray}
\noindent where $C(n,\:\mathcal{E})$, $C(E,\:\mathcal{E})$ are
positive constants depending only on $n$, $\mathcal{E}$ and on
$E$, $\mathcal{E}$ respectively. Let us recall
\begin{eqnarray}\label{stima differenza A(-1)}
&
&\big(A(x^{0},\:a(x^{0}))^{-1}-A(x^{0},\:b(x^{0}))^{-1}\big)\:(x-z_{\tau})\cdot(x-z_{\tau})\nonumber\\
&=&\Bigg(\int_{a(x^{0})}^{b(x^{0})} D_{t}(x^{0},t)^{-1}\:dt\Bigg)\:(x-z_{\tau})\cdot(x-z_{\tau})\nonumber\\
&=&\Bigg(\int_{a(x^{0})}^{b(x^{0})} -A(x^{0},\:t)^{-1}D_{t}(x^{0},t)A(x^{0},\:t)^{-1}\:dt\Bigg)\:(x-z_{\tau})\cdot(x-z_{\tau})\nonumber\\
&\geq& \int_{a(x^{0})}^{b(x^{0})}\mathcal{F}^{-2}\lambda^{-2}\vert
x-z_{\tau}\vert^{2}\:dt,
\end{eqnarray}
\noindent where the ellipticity and the monotonicity assumptions
\eqref{ellitticita'}, \eqref{monotonia} had been used to obtain
the lower bound estimate in \eqref{stima differenza A(-1)}. By
recalling \eqref{stima H1} and combining \eqref{stima A(-1)} with
\eqref{stima differenza A(-1)}, we finally obtain
\begin{eqnarray*}
C(n,\:\mathcal{E})\mathcal{F}^{-2}\lambda^{-2}\parallel
a-b\parallel_{L^{\infty}(\bar\Gamma_{\rho})}\tau^{2-n} &\leq&
C(E,\:\mathcal{E})\big\{\tau^{2-n+\beta}+\tau^{2-n+\beta}+C_{1}\big\}\nonumber\\
&+&
C_{2}\parallel\Lambda_{A(x,\:a)}^{\Gamma}-\Lambda_{A(x,\:b)}^{\Gamma}\parallel_{*}\tau^{2-n},
\end{eqnarray*}
\noindent where $C_{2}$ is a positive constant depending only on
$\textnormal{diam}({\Omega})$, $\lambda$, L, r, h and $\tau^{0}$.
Consequently
\begin{eqnarray*}
\parallel a-b\parallel_{L^{\infty}(\bar\Gamma_{\rho})}\leq C_{2}\:f(\tau)
+C_{3}\parallel\Lambda_{A(x,\:a)}^{\Gamma}-\Lambda_{A(x,\:b)}^{\Gamma}\parallel_{*},
\end{eqnarray*}
\noindent where $C_{2}>0$ is a constant depending only on
$n,\:\lambda,\:E,\:\mathcal{E}$ and $\mathcal{F}$, $C_{3}>0$ is a
constant depending only on
$n,\:\lambda,\:\mathcal{E},\:\mathcal{F},\:\textnormal{diam}({\Omega}),\:L,\:
r,\:h$ and $\tau^{0}$ and $f(\tau)\rightarrow 0$ as
$\tau\rightarrow 0$. If we let $\tau\rightarrow 0$ we obtain
\eqref{stabilita' anisotropa}. $\quad\blacksquare$

We shall need three technical lemmas before we proceed with the
proof of Theorem \ref{stabilita' puntuale}.\\
Given $a\in L^{\infty}(\Omega)$ satisfying the ellipticity
condition \eqref{ellitticita'sigma} and such that it is continuous
in $\overline U_{\rho}^{i}$ we can extend $a$ to all of $\mathbb
R^{n}$ in such a way that the ellipticity conditions are preserved
and $a$ is uniformly continuous in $U_{\rho}$. We shall continue
to call $a$ such an extended function. Let us denote by $\omega$
the modulus of continuity of $a$ in $U_{\rho}$ that is
\begin{equation}\label{modulo di continuita'}
\vert a(x)-a(y)\vert\leq\omega(\vert x-y\vert),
\qquad\textnormal{for}\:\textnormal{any}\quad
x,\:y\in\overline{U}_{\rho},
\end{equation}
\noindent where $\omega$ is a non negative real-valued function on
$\mathbb R^{+}$ so that $\omega(t)\rightarrow 0$ as $t\rightarrow
0^{+}$. Let $\phi_{\varepsilon}$, $\varepsilon>0$, be a usual
family of mollifying kernels with
$\textnormal{supp}\phi_{\varepsilon}\subset B_{\varepsilon}(0)$.
We introduce the mollification of $a$ as
\begin{equation*}
a_{\varepsilon}=\phi_{\varepsilon}\ast a(x).
\end{equation*}
%
%
\begin{lemma}\label{lemma stima con mod. continuita'}
For any $\varepsilon\leq\rho/2$ we have
\begin{equation}\label{stima con modulo di continuita'}
\vert a_{\varepsilon}(x)-a(x)\vert\leq\omega(\varepsilon),
\qquad\textnormal{for}\:\textnormal{any}\quad x\in U_{\rho/2},
\end{equation}
where $a_{\varepsilon}$ is the mollified function of step
$\varepsilon>0$.
\end{lemma}
%
%
\textit{Proof of Lemma \ref{lemma stima con mod. continuita'}.} We
have, for every $x\in\mathbb R^{n}$,
\begin{equation}
a_{\varepsilon}(x)-a(x)=\int_{\vert y-x\vert\leq\varepsilon}
\phi_{\varepsilon}(x-y)\:\big(a (y)-a(x) \big)\:dy,
\end{equation}
and when $x\in U_{\rho/2}$, $\varepsilon\leq\rho/2$, $\vert
y-x\vert\leq\varepsilon$ implies $y\in U_{\rho}$, hence
\begin{eqnarray*}
\vert a_{\varepsilon}(x)-a(x)\vert &\leq&\int_{\vert
y-x\vert\leq\varepsilon}
\vert\phi_{\varepsilon}(x-y)\vert\:\omega(\varepsilon)\:dy\nonumber\\
&=&\omega(\varepsilon),
\qquad\textnormal{for}\:\textnormal{any}\quad x\in U_{\rho/2}.
\quad\blacksquare
\end{eqnarray*}
%
%
\begin{lemma}\label{lemma conv. DN locali}
With the same assumptions as above
\begin{equation}\label{convergenza mappe DN locali}
\parallel\Lambda_{A(\cdot,a_{\varepsilon}(\cdot))}^{\Gamma}-\Lambda_{A(\cdot,a(\cdot))}^{\Gamma}\parallel_{\ast}
\longrightarrow 0,\qquad\textnormal{as}\quad\varepsilon\rightarrow
0^{+},
\end{equation}
\end{lemma}
%
%
\textit{Proof of Lemma \ref{lemma conv. DN locali}.} Let $\phi\in
H^{\frac{1}{2}}_{co}(\Gamma)$, $0<\varepsilon\leq\rho/2$ and take
$u,\:u_{\varepsilon}\in H^{1}(\Omega)$ solutions to the problems
\begin{displaymath}
\left\{ \begin{array}{ll}
\mbox{div}(A(x,a)\nabla{u})=0 & \textrm{$\textnormal{in}\quad\Omega$}\\
u=\phi & \textrm{$\textnormal{on}\quad{\partial\Omega}$}
\end{array} \right.
\end{displaymath}
\noindent and
\begin{displaymath}
\left\{ \begin{array}{ll}
\mbox{div}(A(x,a_{\varepsilon})\nabla{u_{\varepsilon}})=0 &
\textrm{$\textnormal{in}\quad\Omega$}\\
u_{\varepsilon}=\phi &
\textrm{$\textnormal{on}\quad{\partial\Omega}$}
\end{array} \right.
\end{displaymath}
\noindent respectively, then by \eqref{def DN locale} (see
\cite[(b), p. 253]{A}) we have
\begin{eqnarray}\label{Aless. id. lemma 4.2}
\langle \big(\Lambda_{A(\cdot,a_{\varepsilon}(\cdot))}^{\Gamma}-
\Lambda_{A(\cdot,a(\cdot))}^{\Gamma}\big)\phi,\:\phi\rangle
&=&\!\!\!\!\int_{\Omega}(A(x,a^{\varepsilon})-A(x,a))\nabla u_{\varepsilon}\cdot\nabla u\nonumber\\
&=&\!\!\!\!\int_{U_{\rho/2}^{i}}(A(x,a_{\varepsilon})-A(x,a)\nabla u_{\varepsilon}\cdot\nabla u\nonumber\\
&+&\!\!\!\!\!\!\int_{\Omega\setminus
U_{\rho/2}^{i}}\!\!\!\!\!\!(A(x,a_{\varepsilon})\!\!-\!A(x,a))\nabla
u_{\varepsilon}\!\!\cdot\!\nabla u
\end{eqnarray}
\noindent and by combining the H\"older continuity of $A(x,t)$
with \textit{Lemma \ref{lemma stima con mod. continuita'}}
\begin{eqnarray}\label{stima int su U}
\int_{U_{\rho/2}^{i}}(a_{\varepsilon}-a)\nabla
u_{\varepsilon}\cdot\nabla u &\leq&
C(\mathcal{F})\omega(\varepsilon)\parallel\nabla
u_{\varepsilon}\parallel_{L^{2}(\Omega)}
\parallel\nabla u\parallel_{L^{2}(\Omega)}
\nonumber\\
&\leq& \tilde
C\:C(\mathcal{F})\omega(\varepsilon)\parallel\phi\parallel_{H_{co}^{\frac{1}{2}}(\Gamma)}^{2},
\end{eqnarray}
where $C(\mathcal{F})$ is a positive constant depending on the
constant of regularity $\mathcal{F}$ for $A(x,t)$ and $\tilde C$
is a positive constant which does not depend on $\varepsilon$. For
any real numbers $p$, $q$ with $\frac{1}{p}+\frac{2}{q}=1$
\begin{eqnarray}\label{stima int fuori da U}
\!\!\!\!\int_{\Omega\setminus
U_{\rho/2}^{i}}\!\!\!\!\!\!(A(x,a_{\varepsilon}(x)-A(x,a(x)))\nabla
u_{\varepsilon}\cdot\nabla u\leq\!C(\mathcal{F})\!\parallel
a_{\varepsilon}\!-a\parallel_{L^{p}(\Omega)}& &\nonumber\\
\cdot\parallel\nabla
u_{\varepsilon}\parallel_{L^{q}(\Omega\setminus
U^{i}_{\rho/2})}\:\parallel\nabla
u\parallel_{L^{q}(\Omega\setminus U^{i}_{\rho/2})}& &.
\end{eqnarray}
\noindent By Meyers' inequality [M] we have that there exists
$q>2$ such that
\begin{equation}\label{Meyer}
\parallel\nabla u\parallel_{L^{q}(\Omega\setminus U_{\rho})}\leq\:C\parallel\nabla u\parallel_{L^{2}(\Omega)}
\end{equation}
\noindent and the same holds for $u_{\varepsilon}$ and combining
\eqref{stima int su U}-\eqref{Meyer} we obtain
\begin{eqnarray*}
\parallel\Lambda_{A(\cdot,a(\cdot))}^{\Gamma}-\Lambda_{A(\cdot,a_{\varepsilon}(\cdot))}^{\Gamma}\parallel_{\ast}
&=& \sup _{\phi\in
H_{co}^{\frac{1}{2}}(\Gamma),\:\parallel\phi\parallel_{H_{co}^{\frac{1}{2}}(\Gamma)}=1}
\langle\big(\Lambda_{A(\cdot,a(\cdot))}^{\Gamma}-\Lambda_{A(\cdot,a_{\varepsilon}(\cdot))}^{\Gamma}\big)
\phi,\:\phi\rangle\nonumber\\
&\leq&C\big(\omega(\varepsilon)+\parallel
a_{\varepsilon}-a\parallel_{L^{p}(\Omega)}\big),
\end{eqnarray*}
\noindent where $C$ is a positive constant independent from
$\varepsilon$. The above inequality holds for any $\varepsilon$,
$0<\varepsilon\leq\rho/2$, which concludes the proof.
$\quad\blacksquare$
%
%

\textit{Proof of Theorem \ref{stabilita' puntuale}.} Let
$x\in\Gamma_{\rho}$ and take $0<\varepsilon\leq\rho/2$. We can
split the quantity $\vert A(x,a(x))-A(x,b(x))\vert$ as follows
\begin{eqnarray*}
\vert A(x,a(x))-A(x,b(x))\vert &\leq &\vert
A(x,a(x))-A(x,a_{\varepsilon}(x))\vert\nonumber\\
&+& \vert
A(x,a_{\varepsilon}(x))-A(x,b_{\varepsilon}(x))\vert\nonumber\\
&+& \vert A(x,b_{\varepsilon}(x))-A(x,b(x))\vert
\end{eqnarray*}
\noindent and by the H\"older continuity of $A(x,t)$, Lemma
\ref{lemma stima con mod. continuita'} and Theorem \ref{stabilita'
al bordo}
\begin{equation}
\vert A(x,a(x))-A(x,b(x))\vert\leq
C\big(2\omega(\varepsilon)+\:\parallel\Lambda_{a_{\varepsilon}}^{\Gamma}-
\Lambda_{b_{\varepsilon}}^{\Gamma}\parallel_{\ast}\big),
\end{equation}
\noindent where $C$ is a positive constant which does not depend
on $\varepsilon$. By letting $\varepsilon\rightarrow 0^{+}$ and
Lemma \ref{lemma conv. DN locali}, we obtain the desired
estimate. $\quad\blacksquare$\\
Proofs of Theorems \ref{stabilita' Holder al bordo}, \ref{unicita'
al bordo}, \ref{unicita' all'interno} follow the same line of
proofs of Theorems 2.2, 2.3, 2.4 of \cite{AG} by replacing the
singular solutions of Theorem \ref{teore soluzioni singolari} with
the singular solutions with compact support in $\Gamma$ obtained
in Theorem \ref{teore soluzioni singolari 2}.
\subsection{\normalsize{The N-D map.}}\label{N-to-D proofs}
The proof of Theorem \ref{stabilita' puntuale N} shall be based on
the following construction of singular solutions suited for the
\eqref{N bvp} with local data. The following is well known.
\begin{teore}\label{teore funzione neumann}
Let $\Omega$ and $\Gamma$ be as in Theorem \ref{stabilita' al
bordo}. For any $\tau$, $0<\tau\leq\tau_{0}$, set
$z_{\tau}=x^{0}+\tau\tilde\nu$, for some
$x^{0}\in\bar\Gamma_{\rho}$ and $\rho$, $0<\rho\leq\rho_{0}$. If
$L$ is the operator  of \eqref{operatore L}, with H\"older
continuous coefficients matrix
$\sigma=\{\sigma_{ij}\}_{i,\:j=1\dots n}$, with exponent
$0<\beta<1$, the Neumann's function $N_{\sigma}$ for the boundary
value problem associated to the operator \eqref{operatore L}
\begin{displaymath}
\left\{ \begin{array}{ll}
L\:N_{\sigma}(x,\:z_{\tau})=-\delta(x-z_{\tau}), & \textrm{$\textnormal{in}\quad\Omega_{\rho}$}\\
\sigma\nabla
N_{\sigma}(x,\:z_{\tau})\cdot\nu=\frac{1}{\vert\partial\Omega_{\rho}\vert},
& \textrm{$\textnormal{on}\quad{\partial\Omega_{\rho}}$}
\end{array} \right.
\end{displaymath}
has the form
\begin{equation}\label{funzione di neumann}
N_{\sigma}(x,\:z_{\tau})=
C_{n}\:\big(\det(\sigma(z_{\tau}))\big)^{-1/2}\Big(\sigma^{-1}(z_{\tau})(x-z_{\tau})\cdot(x-z_{\tau})\Big)^{\frac{2-n}{2}}+R(x,\:z_{\tau}),
\end{equation}
\noindent where $C_{n}$ is a suitable dimensional constant and the
remainder $R(x,\:z_{\tau})$ satisfies
\begin{equation}\label{stima su R2}
\vert R(x,\:z_{\tau})\vert+\vert x-z_{\tau}\vert\:\vert
\nabla_{x}R(x,\:z_{\tau})\vert \leq C \vert
x-z_{\tau}\vert^{2-n+\alpha},
\end{equation}
\noindent for every $x\in\Omega_{\rho}$, $\vert
x-z_{\tau}\vert\leq r_{0}$, where $C=C(E)$ is a positive constant
depending on $E$, $r_{0}$ is a positive number which depends only
on the geometry of $\Omega$ and $0<\alpha<\beta$. Moreover
\begin{equation}\label{stima H1}
\|\:N_{\sigma}(\cdot,\:z_{\tau})\|_{H^{1}(\Omega)}\leq
C\tau^{(2-n)/2},\qquad for\:any\quad 0<\tau\leq\tau^{0},
\end{equation}
\noindent where $C>0$ is a constant which only depends on
$diam({\Omega})$, $\lambda$, L, r, h and $\tau^{0}$.
\end{teore}
\textit{Proof}. See the proof of Corollary \ref{corollario sol
sing 1} and \cite[Chapter 1]{Mi}. $\quad\blacksquare$
\begin{teore}\label{soluzioni N zero fuori Gamma}
Let $\Omega$ and $\Gamma$ be as in Theorem \ref{stabilita' al
bordo}. For any $\tau$, $0<\tau\leq\tau_{0}$, set
$z_{\tau}=x^{0}+\tau\tilde\nu$, for some
$x^{0}\in\bar\Gamma_{\rho}$ and $\rho$, $0<\rho\leq\rho_{0}$. If
$\sigma$ is the matrix with entries $\{\sigma_{ij}\}_{i,\:j=1\dots
n}$ in \eqref{operatore L} and $S$ is an open portion of
$\partial\Omega_{\rho}\setminus\partial\Omega$ with positive
distance from $\partial\Omega$, there exists $u\in
H^{1}_{loc}(\bar\Omega_{\rho}\setminus z_{\tau})$ solution to
\begin{displaymath}
\left\{ \begin{array}{lll}
Lu=-\delta(x-z_{\tau}), & \textrm{$\textnormal{in}\quad\Omega_{\rho}$}\\
\sigma\nabla u\cdot\nu=0, &
\textrm{$\textnormal{on}\quad{\partial\Omega_{\rho}}$}\\
\sigma\nabla u\cdot\nu=-\frac{1}{\vert S\vert}, &
\textrm{$\textnormal{on}\quad{S}$.}
\end{array} \right.
\end{displaymath}
Moreover
\begin{equation}\label{stima H1 u}
\|\:u(\cdot,\:z_{\tau})\|_{H^{1}(\Omega)}\leq
C\tau^{(2-n)/2}+B,\qquad for\:any\quad0<\tau\leq\tau^{0},
\end{equation}
\noindent where $C,\:B>0$ are constants which only depend on
$diam({\Omega})$, $\lambda$, L, r, h and $\tau^{0}$.
\end{teore}
\textit{Proof}. Let $N_{\sigma}(\cdot,z_{\tau})$ be the Neumann
function for $\Omega_{\rho}$
\begin{displaymath}
\left\{ \begin{array}{ll}
\mbox{div}\big(\sigma\nabla N_{\sigma}(x,z_{\tau})\big)=-\delta(x-z_{\tau}), & \textrm{$\textnormal{in}\quad\Omega_{\rho}$}\\
\sigma\nabla
N_{\sigma}(x,z_{\tau})\cdot\nu=-\frac{1}{\vert\partial\Omega_{\rho}\vert},
& \textrm{$\textnormal{on}\quad{\partial\Omega_{\rho}}$}
\end{array} \right.
\end{displaymath}
and $S$ be an open portion of
$\partial\Omega_{\rho}\setminus\partial\Omega$ with positive
distance from $\partial\Omega$. Set
\begin{equation*}
u(x)=N_{\sigma}(x,z_{\tau})+w(x),\quad\textnormal{for\:any}\quad
x\in\Omega_{\rho},
\end{equation*}
where $w\in H^{1}(\Omega_{\rho})$ is the solution to
\begin{displaymath}
\left\{ \begin{array}{lll}
\mbox{div}\big(\sigma\nabla w)\big)=0, & \textrm{$\textnormal{in}\quad\Omega_{\rho}$}\\
\sigma\nabla w\cdot\nu=-\frac{1}{\vert\partial\Omega_{\rho}\vert},
&
\textrm{$\textnormal{on}\quad{\partial\Omega_{\rho}}\setminus S$}\\
\sigma\nabla w\cdot\nu=-\frac{\vert\partial\Omega_{\rho}\setminus
S\vert}{\vert S\vert\:\vert \partial\Omega_{\rho}\vert}, &
\textrm{$\textnormal{on}\quad{S}$.}
\end{array} \right.
\end{displaymath}
$u$ is a solution of the given boundary value problem and by
Caccioppoli inequality it also satisfies \eqref{stima H1 u}.
$\quad\blacksquare$\\
\textit{Proof of Theorem \ref{stabilita' puntuale N}}. It suffices
to follow the arguments of the proof of Theorem \ref{stabilita'
puntuale} by simply replacing the appropriate singular solutions.
$\quad\blacksquare$
\begin{osse}\label{remark finale}
The argument introduced in Theorem \ref{soluzioni N zero fuori
Gamma} also enables to construct singular solutions of the type of
those introduced in Theorem \ref{teore soluzioni singolari 2}
which however satisfy the zero Neumann condition on
$\partial\Omega\setminus\Gamma$. By means of such singular
solutions it is rather obvious how the proofs of the remaining
Theorems \ref{stabilita' al bordo}, \ref{unicita' al bordo} and
Corollary \ref{unicita' all'interno} can be adapted when the local
Dirichlet-to-Neumann map is replaced by the local
Neumann-to-Dirichlet.
\end{osse}

\end{document}